# FLATTENING AND SUBANALYTIC SETS IN RIGID ANALYTIC GEOMETRY

T.S. Gardener
Hans Schoutens

ABSTRACT. Let $K$ be an algebraically closed field endowed with a complete non-archimedean norm with valuation ring $R$. Let $f\colon Y \to X$ be a map of $K$-affinoid varieties. In this paper we study the analytic structure of the image $f(Y) \subset X$; such an image is a typical example of a subanalytic set. We show that the subanalytic sets are precisely the **D**-semianalytic sets, where **D** is the truncated division function first introduced by DENEF and VAN DEN DRIES. This result is most conveniently stated as a Quantifier Elimination result for the valuation ring $R$ in an analytic expansion of the language of valued fields.

To prove this we establish a Flattening Theorem for affinoid varieties in the style of HIRONAKA, which allows a reduction to the study of subanalytic sets arising from flat maps, i.e., we show that a map of affinoid varieties can be rendered flat by using only finitely many local blowing ups. The case of a flat map is then dealt with by a small extension of a result of RAYNAUD and GRUSON showing that the image of a flat map of affinoid varieties is open in the Grothendieck topology.

Using Embedded Resolution of Singularities, we derive in the zero characteristic case a Uniformization Theorem for subanalytic sets: a subanalytic set can be rendered semianalytic using only finitely many local blowing ups with smooth centres. As a corollary we obtain that any subanalytic set in the plane $R^2$ is semianalytic.

**Acknowledgment.** The authors want to thank JAN DENEF for bringing to their attention RAYNAUD's Theorem and GABRIEL CARLYLE for his minute proof reading.

## 1. RIGID ANALYTIC FLATTENING

**1.1. Blowing Ups.** Let $X$ be a rigid analytic variety. The protagonists of this section will be the local blowing up maps and their compositions. For the definition and elementary properties of rigid analytic blowing up maps, we refer to [**Sch 5**]. Suffice it to say here that they are characterised by a universal property rendering a coherent sheaf of ideals invertible, similar to the classical case. Any blowing up map is proper and an isomorphism away from the centre. If its centre is nowhere dense, then it is also surjective. A local blowing up $\pi$ of $X$ is a composition of a blowing up map $\pi'\colon \tilde{X} \to U$ and an open immersion $U \hookrightarrow X$. We will always assume that $U$ is affinoid. If $Z$ is the centre of the blowing up $\pi'$ (and hence in particular a closed analytic subvariety of $U$), then we call $Z$ also the centre of $\pi$ and we will say that $\pi$ is the local blowing up of $X$ with locally closed centre $Z$.

Let $f\colon Y \to X$ be a map of rigid analytic varieties and let $\pi\colon \tilde{X} \to U \hookrightarrow X$ be a local blowing up with centre $Z \subset U$. If $\theta'\colon \tilde{Y} \to f^{-1}(U)$ denotes the blowing up of $f^{-1}(U)$ with centre $f^{-1}(Z)$ and $\theta\colon \tilde{Y} \to Y$ is the composition of the latter with the open immersion $f^{-1}(U) \hookrightarrow Y$ (making $\theta$ into a local blowing up of $Y$) then by universality of the blowing up, there exists a unique map $\tilde{f}\colon \tilde{Y} \to \tilde{X}$, making the following diagram commute

(1)
$$\begin{array}{ccc} \tilde{Y} & \xrightarrow{\theta} & Y \\ \tilde{f} \downarrow & & \downarrow f \\ \tilde{X} & \xrightarrow{\pi} & X. \end{array}$$

This unique map $\tilde{f}$ is called the *strict transform* of $f$ under $\pi$ and the above diagram will be referred to as the *diagram of the strict transform*.

---

In general, we will not be able to work with just a single local blowing up, but we will make use of maps which are finite compositions of local blowing up maps. Therefore, if $\pi\colon \tilde{X} \to X$ is the composite map $\psi_1 \circ \cdots \circ \psi_m$, with each $\psi_{i+1}\colon X_{i+1} \to U_i \hookrightarrow X_i$ a local blowing up map with centre $Z_i$, for $i < m$, (with $X = X_0$ and $\tilde{X} = X_m$), then we define recursively $f_i\colon Y_i \to X_i$ as the strict transform of $f_{i-1}$ under $\psi_i$ where $f_0 = f$ and $Y_0 = Y$. The last strict transform $f_m$ is called the *(final) strict transform* of $f$ under $\pi$ and the other strict transforms $f_i$, for $i < m$, will be referred to as the *intermediate* strict transforms.

In the sequel we will make extensive use of maps $\pi$ of the above type and we will adopt the notation introduced above. When we want to emphasize the dependence of all these data on $\pi$, we may add $\pi$ as a subscript. For instance, the strict transform of $f$ under $\pi$ might then be denoted by $\tilde{f}_\pi$, given by the commutative diagram

$$\begin{array}{ccc} \tilde{Y}_\pi & \xrightarrow{\theta_\pi} & Y \\ \tilde{f}_\pi \downarrow & & \downarrow f \\ \tilde{X}_\pi & \xrightarrow{\pi} & X. \end{array}$$

For us, the following three possible properties of a map $\pi$ as above, will be crucial.

(i) The centres $Z_i$ are nowhere dense.
(ii) The intermediate strict transforms are flat over their centre, i.e., the restriction $f_i^{-1}(Z_i) \to Z_i$ is flat, for $i < m$.
(iii) The final strict transform $\tilde{f}$ of $f$ under $\pi$ is flat.

Our Flattening Theorem states that given a map of affinoid varieties $f\colon Y \to X$, we can find finitely many maps $\pi_1, \ldots, \pi_s$ as above with these three properties (i)-(iii), such that the union of their images covers $X$.

The proof of this theorem is to be found below in **(1.4)**. It is based on a local flattening theorem which has been proved by the first author in [**Gar**], using a construction of the rigid analytic Voûte Etoilée, after HIRONAKA, and on work of the second author in [**Sch 7**] on rigid analytic flatificators, also after HIRONAKA. In order to make the construction of the Voûte Etoilée, it seems necessary to add extra points to the rigid analytic variety, following BERKOVICH. It is not the intention of this paper to give a proof here, nor will we make any attempts to explain properly what is a Berkovich space, an étoile, the Voûte Etoilée or a flatificator. Nevertheless, a sketch of proof is included for the reader familiar with these notions. A full proof can be found in [**Gar**, Theorem 4.3]. We do include, however, a fully detailed proof of how to derive the Rigid Analytic Flattening Theorem **(1.4)** from this theorem.

**1.2. Local Flattening of Berkovich Spaces.** *Let $f\colon \mathbb{Y} \to \mathbb{X}$ be a map of $K$-affinoid (Berkovich) spaces with $\mathbb{X}$ reduced. Pick $x \in \mathrm{Im}(f)$ and let $\mathbb{L}$ be a non-empty compact subset of $f^{-1}(x)$. There exists a finite collection $E$ of maps $\pi\colon \mathbb{X}_\pi \to \mathbb{X}$, with each $\mathbb{X}_\pi$ affinoid, such that the following four properties hold, where we put $\mathbb{X}_0 = \mathbb{X}$, $\mathbb{Y}_0 = \mathbb{Y}$ and $f_0 = f$.*

(i) *Each $\pi \in E$ is the composition $\psi_1 \circ \cdots \circ \psi_m$ of finitely many local blowing ups $\psi_{i+1}\colon \mathbb{X}_{i+1} \to \mathbb{U}_i \hookrightarrow \mathbb{X}_i$ with nowhere dense centre $\mathbb{Z}_i \subset \mathbb{U}_i$, for $i < m$. (Here we have suppressed in our notation the dependence of these data on $\pi$).*
(ii) *For each $\pi \in E$, define inductively $f_{i+1}$ as the strict transform of $f_i$ under the local blowing up $\psi_{i+1}$. Then $f_i^{-1}(\mathbb{Z}_i) \to \mathbb{Z}_i$ is flat, for $i < m$.*
(iii) *For each $\pi \in E$, the final strict transform $f_\pi\colon \mathbb{Y}_\pi \to \mathbb{X}_\pi$ of $f$ under the whole map $\pi$ (which is $f_m\colon \mathbb{Y}_m \to \mathbb{X}_m$ according to our enumeration) given by the strict transform diagram*

(1)
$$\begin{array}{ccc} \mathbb{Y}_\pi & \xrightarrow{\theta} & \mathbb{Y} \\ f_\pi \downarrow & & \downarrow f \\ \mathbb{X}_\pi & \xrightarrow{\pi} & \mathbb{X}, \end{array}$$

*is flat at each point of $\mathbb{Y}_\pi$ lying above a point of $\mathbb{L}$.*
(iv) *The union of all the $\pi(\mathbb{X}_\pi)$, for $\pi \in E$, is a neighbourhood of $x$.*

*Sketch of Proof.* Let $\mathcal{E}(\mathbb{X})$ denote the collection of all maps $\pi\colon \mathbb{X}' \to \mathbb{X}$ which are finitely many compositions of (Berkovich) local blowing up maps. One can define a partial order relation on $\mathcal{E}(\mathbb{X})$ by calling $\psi\colon \mathbb{X}'' \to \mathbb{X}$ smaller than $\pi$, if $\psi$ factors as $\pi q$. We denote this by $\psi \leq \pi$. Such a $q$ is then necessarily unique and must belong to $\mathcal{E}(\mathbb{X}')$. Any two maps $\pi_1, \pi_2 \in \mathcal{E}(\mathbb{X})$ admit a unique minimum or *meet* $\pi_3 \in \mathcal{E}(\mathbb{X})$ with respect to this ordering, denoted by $\pi_1 \wedge \pi_2$. This meet $\pi_3$ is just the strict transform of $\pi_2$ under $\pi_1$ (or vice versa).

An étoile $e$ on $X$ is now defined as a maximal filter on the semi-lattice $\mathcal{E}(\mathbb{X})$ subject to the extra condition that for any $\pi \in e$ we can find $\psi \in e$, with $\psi \leq \pi$ (i.e., $\psi = \pi q$) such that the image of $q$ is relatively compact (i.e., its closure is compact). The collection of all étoiles on $\mathbb{X}$ is called the Voûte Etoilée of $\mathbb{X}$ and is denoted by $\mathcal{E}_\mathbb{X}$. This space is topologised by taking for opens the sets of the form $\mathcal{E}_\pi$ given as the collection of all étoiles on $\mathbb{X}$ containing $\pi\colon \mathbb{X}' \to \mathbb{X}$, for some $\pi \in \mathcal{E}(\mathbb{X})$. In fact, $\mathcal{E}_\pi$ is isomorphic with $\mathcal{E}_{\mathbb{X}'}$. The Voûte Etoilée is Hausdorff in this topology. Moreover, any étoile defines a unique point $x \in \mathbb{X}$ that lies in the image of all maps belonging to the étoile. Hence there is a canonical map $p_\mathbb{X}\colon \mathcal{E}_\mathbb{X} \to \mathbb{X}$, which is a continuous surjection. It is a highly non-trivial result that this map is also proper in the sense that the inverse image of a compact is compact.[1]

Now, returning to the proof of the theorem, let $e$ be any étoile on $X$ lying above $x$, i.e., $x = p_\mathbb{X}(e)$. Next one has to introduce the concept of the flatificator of $f$ at $x$. This will be the largest locally closed subset $\mathbb{Z}$ of $\mathbb{X}$ containing $x$, such that $f$ is flat over it (i.e., the restriction $f^{-1}(\mathbb{Z}) \to \mathbb{Z}$ is flat). Such a flatificator always exists and is moreover stable under base change (i.e., if $g\colon \mathbb{X}' \to \mathbb{X}$ is arbitrary, then $g^{-1}(\mathbb{Z})$ is the flatificator of the base change $\mathbb{Y} \times_\mathbb{X} \mathbb{X}' \to \mathbb{X}'$ at $x'$, for any $x'$ in the fibre above $x$). Let $\psi_1\colon \mathbb{X}_1 \to \mathbb{X}$ be the local blowing up with centre this flatificator $\mathbb{Z}$ and let $f_1\colon \mathbb{Y}_1 \to \mathbb{X}_1$ denote the strict transform of $f$ under $\psi_1$ giving rise to the commutative strict transform diagram

(2)
$$\begin{array}{ccc} \mathbb{Y}_1 & \xrightarrow{\theta_1} & \mathbb{Y} \\ f_1 \downarrow & & \downarrow f \\ \mathbb{X}_1 & \xrightarrow{\psi_1} & \mathbb{X}. \end{array}$$

The fibre $f_1^{-1}(x_1)$ is naturally a closed subvariety of the original fibre $f^{-1}(x)$ after extending the scalars, for any point $x_1 \in \mathbb{X}_1$ lying above $x$. Moreover, in [**Sch 7**] it is shown that, if $f$ was not flat in some point $y \in \mathbb{L}$, then the closed immersion $f_1^{-1}(x_1) \hookrightarrow f^{-1}(x)$ is strict in some point lying above $y$. This property is not violated if one throws out some dense irreducible components of $\mathbb{Z}$ and hence we may assume that $\mathbb{Z}$ is nowhere dense. It then follows that $\psi_1 \in e$. Under the identification $\mathcal{E}_{\psi_1}$ with $\mathcal{E}_{\mathbb{X}_1}$, the étoile $e$ corresponds to an étoile $e_1$ on $\mathbb{X}_1$. This in turn defines a unique point $x_1 = p_{\mathbb{X}_1}(e_1)$ of $\mathbb{X}_1$ and we can repeat the above process by setting $\mathbb{L}_1 = \mathbb{Y}_1 \cap (\{x_1\} \times \mathbb{L})$.

However, this process cannot go on indefinitely, since any decreasing chain of closed subvarieties (the $f_i^{-1}(x_i) \cap \mathbb{L}$) in a compact set ($\mathbb{L}$) must become stationary. So we must come to a point, say after $m$ steps, where the strict transform $f_m$ is flat in every point of $\mathbb{L}_m$. By an application of the compactness of $\mathbb{L}$ and the fact that flatness is open in the source, one can, after possibly taking one more local blowing up (and in fact, one only needs an open immersion) reduce to the case that $f_m$ is flat at each point of $\mathbb{Y}_m$ lying above a point of $\mathbb{L}$. Let $\pi_e$ be the composition $\psi_1 \circ \cdots \circ \psi_m$, so that $\pi_e$ satisfies conditions (**i**)-(**iii**).

To finish the proof, one now uses the fact that $p_\mathbb{X}\colon \mathcal{E}_\mathbb{X} \to \mathbb{X}$ is surjective and proper, in order to find finitely many étoiles $e$ as above and hence finitely many maps $\pi_e$ for which the properties (**i**)-(**iii**) hold, such that the images of these $\pi_e$ cover a compact neighbourhood of $x$, as required.

A more detailed proof is in [**Gar**, Theorem 4.3] . ∎

**1.3. Lemma.** *Let $f\colon Y \to X$ be a continuous map of topological compact Hausdorff spaces. Let $F$ be a closed subset of $X$ and $V$ an open subset of $Y$, such that $f^{-1}(F) \subset V$. Then there exists an open $U$ of $X$ containing $F$, such that $f^{-1}(U) \subset V$.*

---

[1] The reader familiar with Boolean algebra will note the analogy with the Stone space of ultrafilters on a Boolean lattice. Note that no étoile can be a principal filter. The role of the Frechet filter of cofinite sets is played here by the filter consisting of all maps which are compositions of finitely many blowing up maps with nowhere dense centre, namely, each such map is contained in every étoile.

*Proof.* Since $Y \setminus V$ is closed whence compact, so is its image $f(Y \setminus V)$. By assumption the latter is disjoint from $F$, so that (by normality) we can find an open $U$ of $X$ containing $F$ with $U \cap f(Y \setminus V) = \emptyset$. This $U$ is now as required. ∎

After this auxiliary Lemma, we can now derive the Flattening Theorem from **(1.2)**.

**1.4. Flattening Theorem.** *Let $f\colon Y \to X$ be a map of affinoid varieties with $X$ reduced. Then there exists a finite collection $E$ of maps $\pi\colon X_\pi \to X$, with each $X_\pi$ again affinoid such that the following properties hold.*

  (i) *Each $\pi \in E$ is the composition $\psi_1 \circ \cdots \circ \psi_m$ of finitely many local blowing up maps $\psi_{i+1}$ with locally closed nowhere dense centre $Z_i$, for $i < m$.*
  (ii) *For each $\pi \in E$, let $f_i$ denote the strict transform of $f$ after the $i$-th blowing up $\psi_i$, then $f_i^{-1}(Z_i) \to Z_i$ is flat, for $i < m$. The diagram of strict transform is*

$$\begin{array}{ccc} Y_{i+1} & \xrightarrow{\zeta_{i+1}} & Y_i \\ f_{i+1} \downarrow & & \downarrow f_i \\ X_{i+1} & \xrightarrow{\psi_{i+1}} & X_i. \end{array}$$

  (iii) *The strict transform $f_\pi\colon Y_\pi \to X_\pi$ of $f$ under the whole map $\pi$ (which is $f_m$ according to our enumeration) is flat. The diagram of strict transform is*

(1)
$$\begin{array}{ccc} Y_\pi & \xrightarrow{\theta} & Y \\ f_\pi \downarrow & & \downarrow f \\ X_\pi & \xrightarrow{\pi} & X. \end{array}$$

  (iv) *The union of all the $\mathrm{Im}(\pi)$, for $\pi \in E$, equals $X$.*

*Proof.* Let $\mathbb{X} = \mathbb{M}(X)$ and $\mathbb{Y} = \mathbb{M}(Y)$ be the corresponding Berkovich spaces and let us continue to write $f$ for the corresponding map $\mathbb{Y} \to \mathbb{X}$. Fix an analytic point $x$ of $X$, i.e., a point of $\mathbb{X}$. If $x \notin \mathrm{Im}(f)$, we take the blowing up with empty centre, which amounts to taking the identity map. Otherwise, let $\mathbb{L} = f^{-1}(x)$, which is closed in $\mathbb{Y}$ whence compact since $\mathbb{Y}$ is. By **(1.2)**, we can find a finite collection $E_x$ of maps $\pi\colon \mathbb{X}_\pi \to \mathbb{X}$ with $\mathbb{X}_\pi$ affinoid, such that the conditions **(i)**-**(iv)** hold. For each $\pi \in E_x$, let

$$\begin{array}{ccc} \mathbb{Y}_\pi & \xrightarrow{\theta} & \mathbb{Y} \\ f_\pi \downarrow & & \downarrow f \\ \mathbb{X}_\pi & \xrightarrow{\pi} & \mathbb{X} \end{array}$$

be the corresponding strict transform diagram.

By **(iii)** of **(1.2)** we have that the strict transform $f_\pi$ is flat in each point of $\theta^{-1}(f^{-1}(x)) = f_\pi^{-1}(\pi^{-1}(x))$. Let us first show that we can modify the data in such way that $f_\pi$ becomes flat everywhere. Since flatness is open in the source by [**Sch 7**, Theorem 3.8], we can find an open neighbourhood $\mathbb{V}'$ of $f_\pi^{-1}(\pi^{-1}(x))$ in $\mathbb{Y}_\pi$ over which $f_\pi$ is flat. Applying **(1.3)**, we can find an open neighbourhood $\mathbb{U}'$ of $\pi^{-1}(x)$, such that $f_\pi^{-1}(\mathbb{U}') \subset \mathbb{V}'$ and we can find an open neighbourhood $\mathbb{U}$ of $x$ in $\mathbb{X}$, such that $\pi^{-1}(\mathbb{U}) \subset \mathbb{U}'$. The neighbourhood $\mathbb{U}$ can be taken inside the union of all the $\mathrm{Im}(\pi)$, for all $\pi \in E_x$. Set $\mathbb{U}_\pi = \pi^{-1}(\mathbb{U})$. Note that $\mathbb{U}_\pi \hookrightarrow \mathbb{X}_\pi$ is the strict transform of the open immersion $\mathbb{U} \hookrightarrow \mathbb{X}$ under $\pi$. Let $\psi$ be the restriction of $\pi$ to $\mathbb{U}_\pi$. The strict transform of $f$ under $\psi$ is the map

$$f_\pi^{-1}(\mathbb{U}_\pi) \to \mathbb{U}_\pi,$$

which by construction is flat, since

$$f_\pi^{-1}(\mathbb{U}_\pi) \subset f_\pi^{-1}(\mathbb{U}') \subset \mathbb{V}'.$$

This establishes our claim upon replacing $\pi$ by $\psi$.

Hence we may assume that $f_\pi$ is flat. Note also that in the above process, we have not violated condition (**iv**) of **(1.2)**, so that the $\pi(\mathbb{X}_\pi)$, for all $\pi \in E_x$, form a covering of an affinoid neighbourhood $\mathbb{W}_x$ of $x$ in $\mathbb{X}$. We can translate all these diagrams to the rigid analytic setup and assume that the same diagrams hold with the spaces now rigid analytic varieties (see Remark 1 below), where we keep the same names for our spaces and maps, but just replace any blackboard letter, such as $\mathbb{X},\ldots$, by its corresponding roman equivalent $X,\ldots$, denoting the corresponding rigid analytic variety. In particular, (i)-(iii) hold and we show how to obtain (iv).

Let us now vary the analytic point $x$, so that the $W_x$ cover all analytic points of $X$. Therefore already finitely many do so by [**Ber 2**, Lemma 1.6.2]. In particular, there is a finite collection $S$ of analytic points, such that the union of all $\operatorname{Im}(\pi)$, for all $\pi \in E_x$ and all $x \in S$, cover $X$, i.e., condition (iv) is now verified as well. ∎

*Remark 1.* In this translation process from Berkovich data to rigid analytic data, one needs the following. Let $X$ be an affinoid variety and let $\mathbb{X} = \mathbb{M}(X)$ be the corresponding affinoid Berkovich space. Suppose $\pi\colon \tilde{\mathbb{X}} \to \mathbb{U} \hookrightarrow \mathbb{X}$ is the local blowing up with centre $\mathbb{Z}$, where the latter is a closed subspace of the open $\mathbb{U}$. We can find a wide affinoid $V$ of $X$, such that its closure $\mathbb{M}(V)$ in $\mathbb{X}$ is contained inside $\mathbb{U}$. Hence there exists a closed analytic subvariety $Z$ of $V$, such that $\mathbb{M}(Z) = \mathbb{Z} \cap \mathbb{M}(V)$. Let $p\colon \tilde{X} \to V$ be the blowing up of $V$ with this centre $Z$, then $\mathbb{M}(\tilde{X}) \subset \tilde{\mathbb{X}}$ (see [**Gar**, Lemma 2.2] for the details). So in our translation we will replace $\pi$ by the (rigid analytic) local blowing up $\tilde{X} \to V \hookrightarrow X$. Moreover, if $\mathbb{W}$ is an open inside $\mathbb{U}$ such that its closure $\bar{\mathbb{W}}$ is still contained in $\mathbb{U}$, then we can take $V$ such that $\mathbb{W} \subset \mathbb{M}(V)$ and hence

$$\pi^{-1}(\mathbb{W}) \subset \mathbb{M}(\tilde{X}) \subset \tilde{\mathbb{X}}.$$

Note that the local blowing up $\tilde{\mathbb{W}} \to \mathbb{W} \hookrightarrow \mathbb{X}$ of $\mathbb{X}$ with centre $\mathbb{Z} \cap \mathbb{W}$ coincides with the restriction $\pi^{-1}(\mathbb{W}) \to \mathbb{X}$, so that the rigid analytic local blowing up $\tilde{X} \to X$ is sandwiched by the Berkovich local blowing ups $\pi^{-1}(\mathbb{W}) \to \mathbb{X}$ and $\tilde{\mathbb{X}} \to \mathbb{X}$. The picture is

$$\begin{array}{ccccc}
\tilde{\mathbb{W}} & \longrightarrow & \mathbb{W} & \longrightarrow & \mathbb{X} \\
\downarrow & & \downarrow & & \| \\
\mathbb{M}(\tilde{X}) & \longrightarrow & \mathbb{M}(V) & \longrightarrow & \mathbb{X} \\
\downarrow & & \downarrow & & \| \\
\tilde{\mathbb{X}} & \longrightarrow & \mathbb{U} & \longrightarrow & \mathbb{X},
\end{array}$$

where the composite vertical maps are open immersions and the outer composite horizontal maps are local blowing ups.

Moreover, in this way we can maintain in the rigid analytic version all covering properties which were already satisfied in the Berkovich version.

*Remark 2.* Note that we proved something stronger than condition (iv), namely the union of the images of all $\pi \in E$ covers not only all geometric points of $X$, but also all analytic points.

## 2. Subanalytic Sets

**2.1. Definition.** We now introduce the notion of semianalytic and subanalytic sets in rigid analytic geometry. There are essentially two different ways of viewing these objects, one is geometrical in nature and the other is model-theoretic. We give both point of views and leave it to the reader to pick his favourite. In what follows, let $X = \operatorname{Sp} A$ be a reduced affinoid variety (i.e., $A$ has no non-trivial nilpotent elements).

*2.1.1. The Geometric Point of View.*

A subset $\Sigma$ of $X$ is called *globally (rigid) semianalytic* in $X$, if $\Sigma$ is the union of finitely many *basic* subsets, where the latter are of the form

(1) $\quad \{\, x \in X \mid |p_i(x)| \leq |q_i(x)|, \text{ for } i < n \text{ and } |p_i(x)| < |q_i(x)|, \text{ for } n \leq i < m \,\},$

with the $p_i, q_i \in A$. The set $\Sigma$ is just called *(rigid) semianalytic* in $X$, if there exists a finite admissible affinoid covering $\{X_j\}_{j<t}$ of $X$, such that $\Sigma \cap X_j$ is globally semianalytic in $X_j$, for each $j < t$.

The set $\Sigma$ is called *(rigid) subanalytic* in $X$, if there exists a globally semianalytic subset $\Omega$ of $X \times R^N$, for some $N$, such that $\Sigma = \pi(\Omega)$, where $\pi \colon X \times R^N \to X$ is the projection on the first factor. Whereas the collection of all (globally) semianalytic subsets of $X$ is easily seen to be a Boolean algebra, this is no longer obvious at all for the class of subanalytic sets. Recently, LIPSHITZ and ROBINSON gave a proof of this result in [**LR 2**, Corollary 1.6] . Below, we give a short review of their results, since we will make use of them in the proof of our Quantifier Elimination **(2.7)**.

The reader might wonder whether one should not introduce more local versions of subanalyticity, for instance, what about the projection of a semianalytic set which is not globally semianalytic? It follows however quite easily that we would not enlarge the class of sets at all. We extend the notion of semianalytic and subanalytic to an arbitrary quasi-compact rigid analytic variety $X$ as follows: let $\Sigma \subset X$ then $\Sigma$ is semianalytic (respectively, subanalytic) in $X$, if there exists a finite admissible affinoid covering $\{X_i\}_{i<s}$ of $X$, such that each $\Sigma \cap X_i$ is semianalytic (respectively, subanalytic) in $X_i$.

In order to give a neat description of a subanalytic set, it is convenient to introduce a special function $\mathbf{D}$, first introduced by DENEF and VAN DEN DRIES in their paper [**DvdD**] , in which they describe $p$-adic subanalytic sets. Put

$$\mathbf{D} \colon R^2 \to R \colon (a,b) \mapsto \begin{cases} a/b & \text{if } |a| \leq |b| \neq 0 \\ 0 & \text{otherwise.} \end{cases}$$

We define the *algebra $A^{\mathbf{D}}$ of $\mathbf{D}$-functions* on $X$, as the smallest $K$-algebra of $K$-valued functions on $X$ containing $A$ and closed under the following two operations.

  (i) If $p, q \in A^{\mathbf{D}}$, then also $\mathbf{D}(p, q) \in A^{\mathbf{D}}$.
 (ii) If $p \in A\langle T_1, \ldots, T_N \rangle$ and $q_i \in A^{\mathbf{D}}$ with $|q_i| \leq 1$, for $i = 1, \ldots, N$, then also $p(q_1, \ldots, q_N) \in A^{\mathbf{D}}$.

Here, the function $\mathbf{D}(p, q)$ is to be viewed as a pointwise division, i.e., defined by $x \mapsto \mathbf{D}(p(x), q(x))$. Note also that if $p \in A^{\mathbf{D}}$ then $p$ defines a bounded function on $X$ and hence it makes sense to define $|p| = \sup_{x \in X} |p(x)|$. If we allow in the definition of (globally) semianalytic sets also $\mathbf{D}$-functions rather than just elements of $A$, we may now formulate the definition of *(globally) $\mathbf{D}$-semianalytic* sets: the functions appearing in (1) may be elements of $A^{\mathbf{D}}$. The class of globally $\mathbf{D}$-semianalytic sets coincides with the class of $\mathbf{D}$-semianalytic sets. Our main result now will be that a set is $\mathbf{D}$-semianalytic, if and only if, it is subanalytic.

2.1.2. *The Model-Theoretic Point of View.*

If one wants to initiate the model-theoretic study of the field $K$ with its analytic structure, it is more convenient to consider the valuation ring $R$ instead. This is because strictly convergent power series (in $N$ variables) only converge on the unit disk $R^N$. We propose the following language.

The *analytic language* $\mathcal{L}_{\mathrm{an}}$ for $R$ consists of two 2-ary relation symbols $\mathbf{P}_\leq$ and $\mathbf{P}_<$ and an $n$-ary function symbol $F_f$, for every strictly convergent power series $f$ in $n$-variables of norm at most one, i.e., for every $f \in R\langle X_1, \ldots, X_n \rangle$, where $n = 0, 1, \ldots$. The interpretation of $R$ as an $\mathcal{L}_{\mathrm{an}}$-structure is as follows. Each $n$-ary function symbol $F_f$ is interpreted as the corresponding function $f \colon R^n \to R$, defined by the strictly convergent power series $f$ (note that $|f| \leq 1$, so that $f$ is indeed $R$-valued). The relation symbol $\mathbf{P}_\leq$ interprets the subset $\{(x, y) \in R^2 \mid |x| \leq |y|\}$ of $R^2$, and likewise, $\mathbf{P}_<$ describes the subset $\{(x, y) \in R^2 \mid |x| < |y|\}$. Hence, the atomic formulae in this language (or rather, their interpretation in $R$) are of the following three types

(1) $$f(x) = g(x),$$
(2) $$|f(x)| \leq |g(x)|,$$
(3) $$|f(x)| < |g(x)|.$$

Note that the first type can be rewritten as $|f(x) - g(x)| \leq 0$, so that we actually only have to deal with types (2) and (3). One can of course define $\mathbf{P}_<(x, y)$ as $\neg \mathbf{P}_\leq(y, x)$, but the advantage of not doing so is that all formulae can now be made equivalent with positive ones, i.e., without using the negation symbol. One cannot expect $R$ to have elimination of quantifiers in this language, as it has neither in the real or the $p$-adic case (basically the same counterexample, in essence due to OSGOOD, can be used in all three cases, see our Appendix below).

In an attempt to remedy this, we introduce an expansion $\mathcal{L}_{\text{an}}^{\mathbf{D}}$ of $\mathcal{L}_{\text{an}}$ with one new 2-ary function symbol $\mathbf{D}$, which we will interpret in our structure as the function $\mathbf{D}$ of above. If $K$ were the $p$-adic field (and hence $R = \mathbb{Z}_p$), then by a theorem of DENEF and VAN DEN DRIES [**DvdD**], $R$ admits Elimination of Quantifiers in an expansion of this language where one needs to add extra predicates, one for each $n = 2, 3, \ldots$, to express that an element is an $n$-th power; a similar expansion occurs in MACINTYRE's algebraic Quantifier Elimination for $\mathbb{Z}_p$. In the algebraically closed case these predicates are clearly obsolete. Hence the following is the natural rigid analytic analogue: the valuation ring $R$ of $K$ admits Elimination of Quantifiers in the language $\mathcal{L}_{\text{an}}^{\mathbf{D}}$.

Let us see how this ties in with the above notion of subanalyticity. A subset of $R^N$ which is definable in the language $\mathcal{L}_{\text{an}}$ by a quantifier free formula, is precisely a globally semianalytic set whereas an existentially definable set is precisely a subanalytic set. It is not too hard to see that the function $\mathbf{D}$ is existentially definable and whence also every $\mathbf{D}$-function on $R^N$, so that any $\mathbf{D}$-semianalytic subset of $R^N$ is subanalytic. Claiming that the converse also holds is then equivalent with the aforementioned Quantifier Elimination in the language $\mathcal{L}_{\text{an}}^{\mathbf{D}}$.[2]

We remark that any affinoid variety $X$ is quantifier-free definable in $\mathcal{L}_{\text{an}}$ since there is a closed immersion $X \hookrightarrow R^N$ for some $N \in \mathbb{N}$. More generally any quasi-compact rigid analytic variety is also quantifier-free definable in $\mathcal{L}_{\text{an}}$. Also note that semianalytic sets (respectively, subanalytic sets) in such a variety $X$ correspond to quantifier-free definable (respectively, existentially definable) subsets of $X$.

We will be adopting from now on the geometric point of view. An additional advantage is then that we do not really need the field $K$ to be algebraically closed. However, for sake of simplicity we will maintain this assumption in what follows. In particular, one can and we will identify $\text{Sp}(K\langle S_1, \ldots, S_n\rangle)$ with $R^n$.

**2.1.3. Example.** If $f \colon Y \to X$ is a map of affinoid varieties, then the image $f(Y)$ is a typical subanalytic subset of $X$ (not necessarily semianalytic!). Subanalyticity follows from projecting the graph of $f$ (which is analytic, whence semianalytic) onto $X$. More generally, it follows that $f(\Sigma) \subset X$ is subanalytic whenever $\Sigma \subset Y$ is subanalytic. This example shows that even when one is merely interested in closed analytic subsets, one needs to study subanalytic sets as well. However, there are some particular kind of maps which have better understood image. For instance, KIEHL's Proper Mapping Theorem [**Ki**] (or [**BGR**, 9.6.3. Proposition 3]) states that the image of a proper map is closed analytic. However, this does not tell us anything on the image of a semianalytic set under a proper map. In fact, in [**Sch 1**] and [**Sch 2**] the second author shows that if $\Sigma \subset Y$ is semianalytic and $f \colon Y \to X$ is proper, then $f(\Sigma)$ is $\mathbf{D}$-semianalytic in $X$; he carries out a systematic study of the sets arising in this way–*strongly subanalytic* sets. One might hope though that certain proper maps, viz. blowing up maps, nevertheless behave better with respect to semianalyticity. It is the contents of **(2.3)** below that this is true provided one replaces semianalytic by $\mathbf{D}$-semianalytic and then, unfortunately, this is only true generically, i.e., away from the centre. It is because of this (rather straightforward) result that $\mathbf{D}$-functions are needed. Noteworthy here is that in case of the blowing up of the plane in a single (reduced) point, the image of a semianalytic set is nevertheless semianalytic again. This (much harder) result will be used implicitly in the proof of **(3.2)**.

A second class of affinoid maps with well-understood images are the flat maps: their images are finite unions of rational domains and hence in particular semianalytic. This highly non-trivial result is due to RAYNAUD and GRUSON (a full account by MEHLMANN appeared in [**Meh**]). Because of its crucial role in our argument and since we need a slight improvement of their original result in the form **(2.2)** below, we will provide most of the details in Section 4.

**2.2. Theorem (Raynaud-Gruson-Mehlmann).** *Let $f \colon Y \to X$ be a flat map of affinoid varieties. Let $\Sigma$ be a semianalytic subset of $Y$ defined by finitely many inequalities of the form $|h(y)| < 1$ or $|h(y)| \geq 1$, where each $h$ belongs to the affinoid algebra of $Y$ and has supremum norm at most one. Then $f(\Sigma)$ is semianalytic in $X$.*

*Proof.* See Section 4. ∎

---

[2] By an easy logic argument, it is enough to eliminate only existential quantifiers to obtain Quantifier Elimination.

**2.3. Proposition.** *Let $\pi\colon \tilde{X} \to X$ be a map of rigid analytic varieties and let $\Sigma$ be a $\mathbf{D}$-semianalytic subset of $\tilde{X}$. If $\pi$ is a locally closed immersion, then $\pi(\Sigma)$ is $\mathbf{D}$-semianalytic in $X$. If $\pi$ is a local blowing up map with centre $Z$, then $\pi(\Sigma) \setminus Z$ is $\mathbf{D}$-semianalytic in $X$.*

*Proof.* For closed immersions the statement is trivial. If $U = \operatorname{Sp} C \hookrightarrow X = \operatorname{Sp} A$ is a rational affinoid subdomain, then $C = A\langle f/g \rangle$, where $f = (f_1, \ldots, f_n)$ with $f_i, g \in A$ having no common zero. Hence any function $h \in C$ defined on $U$ is $\mathbf{D}$-definable on $X$ (just replace any occurrence of $f_i/g$ by $\mathbf{D}(f_i, g)$). Now, any affinoid subdomain is a finite union of rational subdomains by [**BGR**, 7.3.5. Corollary 3] and hence we proved the proposition for any affinoid open immersion as well. From this the general locally closed immersion case follows easily.

This leaves us with the case of a blowing up. Without loss of generality, we may assume $X$ to be affinoid. Let us briefly recall the construction of a blowing up map as described in [**Sch 5**]. Let $X = \operatorname{Sp} A$ and let $Z$ be a closed analytic subvariety of $X$ defined by the ideal $(g_1, \ldots, g_n)$ of $A$. We can represent $A$ as a quotient of some $K\langle S \rangle$, with $S = (S_1, \ldots, S_m)$, so that $X$ becomes a closed analytic subvariety of $R^m$. However, in order to construct the blowing up of $X$ with centre $Z$, we need a different embedding, given by the surjective algebra morphism

$$K\langle S, T \rangle \twoheadrightarrow A\colon T_j \mapsto g_j,$$

for $j = 1, \ldots, n$, extending the surjection $K\langle S \rangle \twoheadrightarrow A$ and where $T = (T_1, \ldots, T_n)$. This gives us a closed immersion $i\colon X \hookrightarrow R^m \times R^n$ and after identifying $X$ with its image $i(X)$, we see that $Z = X \cap (R^m \times 0)$. Now, the blowing up $\pi\colon \tilde{X} \to X$ is given by a strict transform diagram

$$\begin{array}{ccc} \tilde{X} & \xrightarrow{\pi} & X \\ \tilde{\imath} \downarrow & & \downarrow i \\ W & \xrightarrow{\gamma} & R^m \times R^n \end{array}$$

where $\gamma$ denotes the blowing up of $R^m \times R^n$ with centre the linear space $R^m \times 0$. There is a standard finite admissible affinoid covering $\{W_1, \ldots, W_n\}$ of $W$ where each $W_j$ has affinoid algebra

$$C_j = \frac{K\langle S, T, U \rangle}{(T_j U_1 - T_1, \ldots, T_j U_n - T_n)},$$

so that $\gamma(s, t, u) = (s, t)$ for any point $(s, t, u) \in W_j$, where the latter is considered as a closed analytic subset of $R^m \times R^n \times R^n$ via the above representation of $C_j$. Moreover, $\tilde{X}$ is a closed analytic subvariety of $X \times_{(R^m \times R^n)} W$. Therefore, if we set $\tilde{X}_j = \tilde{\imath}^{-1}(W_j)$, then $\{\tilde{X}_1, \ldots, \tilde{X}_n\}$ is a finite admissible affinoid covering of $\tilde{X}$ with the affinoid algebra $\tilde{A}_j$ of each $\tilde{X}_j$ some quotient of the affinoid algebra

(1) $$\frac{A\langle \hat{U}_j \rangle}{(g_j U_1 - g_1, \ldots, g_j U_n - g_n)}$$

of $W_j \times_{(R^m \times R^n)} X$, where $\hat{U}_j$ means all variables $U_k$ save $U_j$.

With this notation, let us return to the proof of the proposition. We are given some $\mathbf{D}$-semianalytic set $\Sigma$ of $\tilde{X}$ and we seek to describe the image $\pi(\Sigma) \setminus Z$. Let us focus for the time being at one $\Sigma \cap \tilde{X}_j$, where $j \in \{1, \ldots, n\}$. Since $\Sigma \cap \tilde{X}_j$ is $\mathbf{D}$-semianalytic, we can find a quantifier free $\mathcal{L}_{\mathrm{an}}^{\mathbf{D}}$-formula $\varphi(\bar{s}, \bar{u})$, such that $(s, u) \in R^m \times R^n$ belongs to $\Sigma \cap \tilde{X}_j$, if and only if, $\varphi(s, u)$ holds. Hence, for $s \in R^m$, we have that $s \in \pi(\Sigma \cap \tilde{X}_j)$, if and only if,

(2) $$(\exists \bar{u}) \varphi(s, \bar{u}).$$

Note that by (1), if $\varphi(s, u)$ holds, then in particular $(s, u) \in \tilde{X}_j$ and hence $g_j(s) u_k = g_k(s)$, for all $k = 1, \ldots, n$. Now, a point $s \in R^m$ does not belong to $Z$, precisely when one of the $g_k(s)$ does not vanish. Therefore, as $j$ ranges through the set $\{1, \ldots, n\}$ and using (2), it is not too hard to see that $s \in R^m$ belongs to $\pi(\Sigma) \setminus Z$, if and only if, $s \in X$ and

$$\bigvee_{j=1}^{n} \bigwedge_{k=1}^{n} |g_k(s)| \leq |g_j(s)| \wedge g_j(s) \neq 0 \wedge \varphi(s, \mathbf{D}(g_1(s), g_j(s)), \ldots, \mathbf{D}(g_n(s), g_j(s))),$$

which is indeed a **D**-semianalytic description of $\pi(\Sigma) \setminus Z$. ∎

*Remark.* The above result is unsatisfactory in so far as it does not tell us anything about $\pi(\Sigma)$ restricted to the centre $Z$ of the blowing up. If we could prove that also $\pi(\Sigma) \cap Z$ were **D**-semianalytic, then the whole image $\pi(\Sigma)$ would be **D**-semianalytic, as we would very much like to show. But, above $Z$, the map $\pi$ looks like a projection map, so that we can't say much more about $\pi(\Sigma \cap \pi^{-1}(Z)) = \pi(\Sigma) \cap Z$ except that it is a subanalytic set. If $Z$ would be zero dimensional and whence finite, then clearly also $\pi(\Sigma) \cap Z$ is **D**-semianalytic. This suggests that we might be able to use the above result in order to prove Quantifier Elimination by an induction argument on the dimension of $X$, as soon as we can arrange that $Z$ has strictly smaller dimension than $X$. This will be the case, if $Z$ is nowhere dense; a condition we ensure will always be fulfilled.

Another point ought to be mentioned here: although a blowing up $\pi \colon \tilde{X} \to X$ is an isomorphism outside its centre $Z$, this does *not* automatically imply that one can deduce from the **D**-semianalyticity of $\Sigma \setminus \pi^{-1}(Z)$ the same property for its (isomorphic) image $\pi(\Sigma) \setminus Z$. What is going on here is that being (**D**-)semianalytic is not an intrinsic property of a set, but of its embedding in a larger space. In other words, being isomorphic as point sets is not enough and thus the above statement is not a void one.

Before we turn to the proof of our main theorem, let us give a brief review on the model-completeness result of LIPSHITZ and ROBINSON. Geometrically, this amounts to the fact that the complement of a subanalytic set is again subanalytic. This is by no means a straightforward result. In the real case it was shown by GABRIELOV using quite involved arguments and it was only since the appearance of the seminal paper [**DvdD**] of DENEF and VAN DEN DRIES that one has a conceptual proof through a much stronger result, namely, the class of subanalytic sets is equal to the class of **D**-semianalytic sets. Closure under complementation is now immediate. Using the result of [**LR 2**] we exploit their dimension theory to prove our main Quantifier Elimination Theorem (the formulation of which is entirely in the style of [**DvdD**], but the proof of which uses completely different methods, more in the style of HIRONAKA).

**2.4. Theorem (Lipshitz-Robinson).** *The complement $X \setminus \Sigma$ and the closure $\bar{\Sigma}$ (in the canonical topology) of a subanalytic set $\Sigma$ in $X$, where $X$ is a reduced quasi-compact rigid analytic variety, is again subanalytic.*

**2.5. Theorem (Lipshitz-Robinson).** *Let $X$ be a reduced quasi-compact rigid analytic variety and let $\Sigma$ be a subanalytic set in $X$. Then there exists a finite partition of $\Sigma$ by pairwise disjoint rigid analytic submanifolds $X_i$ of $X$ such that their underlying set is subanalytic in $X$.*

The proofs of both Theorems rely on a certain Quantifier Elimination result in some appropriate language and we refer the reader to the paper [**LR 2**, Corollary 1.2 and 1.3] by LIPSHITZ and ROBINSON. Let us just show how one can derive a good dimension theory for subanalytic sets from these results. First of all, there is the notion of the dimension of a quasi-compact rigid analytic variety. This is defined as the maximum of the (Krull) dimension of all its local rings (we give the empty space dimension $-\infty$). In case $X = \operatorname{Sp} A$ is affinoid, this is just the dimension of $A$. Next, we define the dimension of a subanalytic set $\Sigma$ in $X$ as the maximum of all $\dim Y$, where $Y \subset \Sigma$ is a submanifold of $X$. If $\Sigma$ carries already the structure of a manifold, then clearly its subanalytic dimension equals its manifold dimension.

The relevant properties for this dimension function are now summarized by the following proposition.

**2.6. Proposition.** *Let $X$ be a quasi-compact rigid analytic variety and let $\Sigma$ and $\Sigma'$ be (non-empty) subanalytic sets in $X$. Then the following holds.*
  (i) *If $\Sigma \subset \Sigma'$, then the dimension of $\Sigma$ is at most the dimension of $\Sigma'$.*
 (ii) *The dimension of $\Sigma$ is zero, if and only if, $\Sigma$ is finite.*
(iii) *The dimension of $\Sigma$ equals the dimension of its closure (in the canonical topology) $\bar{\Sigma}$.*
(iv) *The dimension of the boundary $\bar{\Sigma} \setminus \Sigma$ is strictly smaller than the dimension of $\Sigma$.*
 (v) *If $f \colon X \to Y$ is a map of quasi-compact rigid analytic varieties, then the dimension of $f(\Sigma)$ is at most the dimension of $\Sigma$, with equality in case $f$ is injective.*
(vi) *If $\Sigma$ is semianalytic, then the dimension of $\Sigma$ is equal to the (usual) dimension of its Zariski closure.*

*Remark.* Note that by **(2.4)** both the closure $\bar{\Sigma}$ and the boundary $\bar{\Sigma} \setminus \Sigma$ are indeed subanalytic.

*Proof.* The first two statements follow from the fact that the dimension of a subanalytic set is the maximum of the dimensions of each manifold in any finite subanalytic manifold partitioning (as in **(2.5)**). The other statements require more work. See [**Lip**] and also [**DvdD**, 3.15-3.26] for the $p$-adic analogues–the proofs just carry over to our present situation, once one has **(2.5)**. ∎

**2.7. Theorem (Quantifier Elimination).** *Let $X$ be a reduced affinoid variety, then the subanalytic subsets of $X$ are precisely the **D**-semianalytic subsets of $X$.*

*Proof.* We have already seen that **D**-semianalytic sets are subanalytic. To prove the converse, let $\Sigma$ be a subanalytic set of $X$. We will induct on the dimension of $\Sigma$ and then on the dimension of $X$. The zero-dimensional case follows immediately from (ii) in **(2.6)**. Hence fix $\dim \Sigma = k > 0$ and $\dim X = d > 0$.

*Step 1.* It suffices to take $\Sigma$ closed in the canonical topology (i.e., the induced topology coming from the norm). Indeed, assume the theorem proven for all subanalytic sets which are closed in the canonical topology. Let $\bar{\Sigma}$ be the closure of $\Sigma$ with respect to the canonical topology. By **(2.4)** and (iii) of **(2.6)**, also $\bar{\Sigma}$ is subanalytic and of dimension equal to the dimension of $\Sigma$. Hence by our assumption $\bar{\Sigma}$ is even **D**-semianalytic. Let $\Gamma$ be the boundary $\bar{\Sigma} \setminus \Sigma$, which is again subanalytic by **(2.4)**. Moreover, by (iv) of **(2.6)**, $\Gamma$ has strictly smaller dimension than $\Sigma$. Hence, by our induction hypothesis on the dimension of a subanalytic set, we have that also $\Gamma$ is **D**-semianalytic. Therefore also $\Sigma = \bar{\Sigma} \setminus \Gamma$, as required.

*Step 2.* Hence we may assume that $\Sigma$ is closed in the canonical topology. There exists a globally semianalytic subset $\Omega' \subset X \times R^N$, for some $N$, such that $\Sigma = f'(\Omega')$, where $f' \colon X \times R^N \to X$ is the projection on the first factor. The union of finitely many **D**-semianalytic sets is again such. Therefore, without loss of generality, we may even take $\Omega'$ to be a basic set, i.e., of the form

$$\left\{ (x,t) \in X \times R^N \mid \bigwedge_{i < m} |p_i(x,t)| \leq |q_i(x,t)| \ \wedge \ \bigwedge_{m \leq i < n} |p_i(x,t)| < |q_i(x,t)| \right\},$$

where the $p_i$ and $q_i$ are in $A\langle T \rangle$, with $X = \mathrm{Sp}\, A$ and $T = (T_1, \ldots, T_N)$. Introduce $n$ new variables $Z_i$ and consider the following closed analytic subset $Y$ of $X \times R^{N+n}$ given by the equations $p_i - Z_i q_i = 0$, for $i < n$. Let $\Omega$ be the basic subset of $Y$ given by $(x, t, z) \in Y$ belongs to $\Omega$ whenever $|z_i| < 1$, for $m \leq i < n$. Let $q$ be the product of all the $q_i$, for $m \leq i < n$, and we obviously can assume that $q \neq 0$ lest $\Sigma$ is non-empty. If $f \colon Y \to X$ denotes the composition of the closed immersion $Y \hookrightarrow X \times R^{N+n}$ followed by the projection $X \times R^{N+n} \to X$, then $f(\Omega \cap U) = \Sigma$, where $U$ is the complement in $Y$ of the zero-set of $q$. Using [**Sch 6**, Corollary 2.2], we may, after perhaps modifying some of the equations defining $Y$, assume that the closure of $U$ in the canonical topology equals the whole of $Y$ and hence the closure (in the canonical topology) of $\Omega \cap U$ is $\Omega$. Now $\Omega \cap U \subset f^{-1}(\Sigma)$ and so $\Omega = \overline{\Omega \cap U} \subset f^{-1}(\Sigma)$, since $\Sigma$ is closed and $f$ is continuous. Hence $f(\Omega) = \Sigma$.

Before giving the details of the proof, let's pause to give a brief outline of how we will go about. According to our Flattening Theorem, we can find finitely many diagrams

$$(\dagger)_\pi \qquad \begin{array}{ccc} Y_\pi & \xrightarrow{\theta_\pi} & Y \\ f_\pi \downarrow & & \downarrow f \\ X_\pi & \xrightarrow{\pi} & X \end{array}$$

indexed by maps $\pi$, where each such $\pi$ is a finite composition of local blowing up maps with the properties (i)-(iii) and such that $X$ is covered by the union of all the $\mathrm{Im}(\pi)$. Now, in order to study $\Sigma = f(\Omega)$, we will chase $\Omega$ around these diagrams $(\dagger)_\pi$. There are only finitely many $\pi$ to consider; it will suffice to do this for one such $\pi$ since the analysis for the others is identical. First we take the preimage $\theta_\pi^{-1}(\Omega)$, which is again a semianalytic set defined by inequalities of the form $|h| < 1$ where the $h$ are functions on $Y_\pi$ of supremum norm at most one. Next we take the image of the latter set under $f_\pi$. Our extension of Raynaud's Theorem **(2.2)** guarantees that this image is semianalytic. Finally we push this set back to $X$ via $\pi$ and denote this set temporarily by $\Sigma'$. If we had the full version of **(2.3)**, i.e., a local blowing up map preserves **D**-semianalyticity, then this last set would be indeed **D**-semianalytic.

Of course, in chasing $\Omega$ around the diagram, we might have lost some points, i.e., it may well be the case that $\Sigma' \neq \Sigma$. But this could happen only for points coming from one of the centres of the local blowing ups that make up $\pi$ (since outside its centre, a blowing up map is an isomorphism). Above each of these centres the strict transform is flat so we account for those missing points using **(2.2)** once more. Hence the only problem in the above reasoning lies in the application of **(2.3)**: it is not the whole image that we can account for by means of that proposition, but only for the part outside the centre. However, the latter has dimension strictly smaller and by an induction argument on the dimension, we could also deal with this part.

*Step 3.* Our second induction hypothesis says that any subanalytic set in an affinoid variety of dimension strictly smaller than $d$ is **D**-semianalytic. Let us first draw the following strengthening of **(2.3)**:

**(2.3)'** *Let $\pi\colon \tilde{W} \to W$ be any local blowing up of a quasi-compact rigid analytic variety $W$ of dimension at most $d$ whose centre $Z$ is nowhere dense. If $\Gamma \subset \tilde{W}$ is **D**-semianalytic, then $\pi(\Gamma) \subset W$ is also **D**-semianalytic.*

The key point is that $Z$ has dimension strictly smaller than the dimension $d$ of $W$, which is also the dimension of $\tilde{W}$. Now
$$\pi(\Gamma) = (\pi(\Gamma) \setminus Z) \cup (Z \cap \pi(\Gamma)).$$
By **(2.3)** we know that $\pi(\Gamma) \setminus Z$ is **D**-semianalytic and by our induction hypothesis on the dimension we have that also $Z \cap \pi(\Gamma)$ is (take a finite affinoid covering to reduce to the affinoid case).

*Step 4.* Now, according to **(1.4)**, there exists a finite collection $E$ of maps $\pi\colon X_\pi \to X$, such that each $\pi \in E$ induces a strict transform diagram $(\dagger)_\pi$ with properties (i)-(iv) of loc. cit. (The intermediate strict transform diagrams are given by $(\dagger)_i$ below). By (iv), if we could show that each $\operatorname{Im}(\pi) \cap \Sigma$ is **D**-semianalytic in $X$, then the same would hold for $\Sigma$, since there are only finitely many $\pi$. Therefore, let us concentrate on one such $\pi = \pi_1 \circ \ldots \circ \pi_m$ and adopt the notation from **(1.1)** for this map, so that in particular, (i)-(iii) of loc. cit. holds. Let each $\pi_{i+1}$ be the blowing up of the admissible affinoid $U_i \subset X_i$ with nowhere dense centre $Z_i \subset U_i$. The diagram of strict transform is given by

$(\dagger)_i$
$$\begin{array}{ccc} Y_{i+1} & \xrightarrow{\theta_{i+1}} & Y_i \\ {\scriptstyle f_{i+1}}\downarrow & & \downarrow{\scriptstyle f_i} \\ X_{i+1} & \xrightarrow[\pi_{i+1}]{} & X_i. \end{array}$$

Define inductively $\Omega_i \subset Y_i$ as $\theta_i^{-1}(\Omega_{i-1})$ starting from $\Omega_0 = \Omega$. Note that each $\Omega_i$ is a semianalytic set of $Y_i$ defined by several inequalities of the type $|h| < 1$, where each $h \in \mathcal{O}(Y_i)$ is of supremum norm at most one. Define also inductively, but this time by downwards induction, the sets $W_{i-1} = \pi_i(W_i) \subset U_i \subset X_i$ where we start with $W_m = X_m = X_\pi$. In particular, we have that $W_0 = \operatorname{Im}(\pi)$. By **(2.3)'** each $W_i$ is **D**-semianalytic in $X_i$. In order to describe $\Sigma$, we will furthermore make use of the sets $\Gamma_i$ defined as $f_i(\Omega_i) \cap W_i$, for $i \leq m$. In particular, note that $\Gamma_0$ is nothing else than $f(\Omega) \cap W_0 = \Sigma \cap \operatorname{Im}(\pi)$, which we aim to show is **D**-semianalytic.

The next claim shows how two successive members in the chain of commutative diagrams $(\dagger)_i$ relate the $\Gamma_i$: for each $i < m$, we have an equality

$(\ddagger)_i$
$$\Gamma_i = \pi_{i+1}(\Gamma_{i+1}) \cup (\Gamma_i \cap Z_i).$$

Assume we have established already $(\ddagger)_i$, for each $i < m$. We will prove, by downwards induction on $i \leq m$, that each $\Gamma_i$ is **D**-semianalytic in $X_i$, so that in particular $\Gamma_0$ would be **D**-semianalytic in $X$, as required. First of all, since $f_\pi = f_m$ is assumed to be flat, we can apply **(2.2)** to $\Omega_m$ to conclude that $\Gamma_m = f_m(\Omega_m)$ is semianalytic whence **D**-semianalytic in $X_m$. Assume now that we have already proven that $\Gamma_{i+1}$ is **D**-semianalytic in $X_{i+1}$ and we want to obtain the same conclusion for $\Gamma_i$ in $X_i$. Using $(\ddagger)_i$, it is enough to establish this for both sets in the right hand side of that equality. The first of these, $\pi_{i+1}(\Gamma_{i+1})$, is **D**-semianalytic since we have now the strong version **(2.3)'** of **(2.3)** at our disposal. As for the second set, $\Gamma_i \cap Z_i$, also this one is **D**-semianalytic, since $f_i$ restricted to $f_i^{-1}(Z_i)$ is flat and since

$$\Gamma_i \cap Z_i = f_i(\Omega_i \cap f_i^{-1}(Z_i)) \cap W_i,$$

so that **(2.2)** applies. Note that we already established that $W_i$ is **D**-semianalytic.

Therefore, it only remains to prove $(\ddagger)_i$. The inclusion $\supset$ is straightforward and we omit the details. To prove $\subset$, let $x_i \in \Gamma_i$. That means that there exists $y_i \in \Omega_i$ and $w_{i+1} \in W_{i+1}$ such that $f_i(y_i) = x_i = \pi_{i+1}(w_{i+1})$. If $x_i \in Z_i$, we are done. Hence assume that $x_i \notin Z_i$ so that $y_i \notin f_i^{-1}(Z_i)$. However, since $W_i \subset U_i$ we have that $y_i \in f_i^{-1}(U_i)$. Since $\theta_{i+1}$ is the blowing up of $f_i^{-1}(U_i)$ with centre $f_i^{-1}(Z_i)$ and whence an isomorphism outside this centre, we can even find $y_{i+1} \in Y_{i+1}$, such that $\theta_{i+1}(y_{i+1}) = y_i$. From $y_i \in \Omega_i$ it then follows that $y_{i+1} \in \Omega_{i+1}$. Put $x_{i+1} = f_{i+1}(y_{i+1})$. Commutativity of the strict transform diagram implies that $\pi_{i+1}(x_{i+1}) = x_i = \pi_{i+1}(w_{i+1})$. Since $x_i \notin Z_i$, the blowing up $\pi_{i+1}$ is an isomorphism in that point, so that $w_{i+1} = x_{i+1}$ which therefore belongs to $f_{i+1}(\Omega_{i+1}) \cap W_{i+1} = \Gamma_{i+1}$, proving our claim, and hence also our main theorem. ∎

*Remark.* We can derive from the above proof also a weak (=non-smooth) uniformization as follows. Define $\Sigma_i$ inductively as the inverse image of $\Sigma_{i-1}$ under $\pi_i$, for $1 \leq i \leq m$, with $\Sigma_0 = \Sigma$. With notations as in the above proof, we can derive, for $i < m$, from $(\ddagger)_i$ the following identity

$$\Sigma_{i+1} \cap W_{i+1} = \Gamma_{i+1} \cup \left(\pi_{i+1}^{-1}(\Gamma_i \cap W_i) \cap W_{i+1}\right).$$

For $i = m - 1$, this takes the simplified form $\Sigma_m = \Gamma_m \cup \pi_m^{-1}(\Gamma_{m-1} \cap Z_{m-1})$. Now, as already observed, $\Gamma_m$ is semianalytic in $X_m = X_\pi$ and similarly $\Gamma_{m-1} \cap Z_{m-1}$ is semianalytic in $X_{m-1}$ and whence also its preimage under $\pi_m$. In other words, we showed the following proposition.

**2.8. Corollary.** *Let $X$ be a reduced affinoid variety and let $\Sigma$ be a subanalytic set in $X$. There exists a finite collection of compositions of finitely many local blowing up maps $\pi_1, \ldots, \pi_n$ with nowhere dense centre, such that the union of the $\operatorname{Im}(\pi_i)$ equals $X$, and such that each preimage $\pi_i^{-1}(\Sigma)$ has become semianalytic.*

*Proof.* This follows from the above discussion in the case where $\Sigma$ is closed in the canonical topology. The reduction to this case uses an induction argument similar to the one in the proof of the theorem. ∎

Note also that to prove the corollary, we do not make use of **(2.3)** but only of **(2.2)**. For an improvement of **(2.8)**, at least in the zero characteristic case, see the Uniformization Theorem **(3.1)** below, where we will be able to take smooth centres for the blowing ups involved.

## 3. Uniformization

In [**Sch 2**, Theorem 4.4] it was proved that for any strongly subanalytic set $\Sigma$ in an affinoid manifold $X$, there exists a finite covering family of compositions $\pi$ of finitely many local blowing ups with smooth and nowhere dense centre, such that the preimage $\pi^{-1}(\Sigma)$ is semianalytic, provided the characteristic of $K$ is zero. The restriction to zero characteristic is entirely due to the lack of an Embedded Resolution of Singularities in positive characteristic. A proof of this rigid analytic Embedded Resolution of Singularities for zero characteristic can be found in [**Sch 4**, Theorem 3.2.5]. In the present paper, we will extend the above Uniformization Theorem to the class of all subanalytic sets. The proof is completely the same as for the strong subanalytic case, in that we only make use of the fact that a subanalytic set is **D**-semianalytic. For the convenience of the reader we give below an outline of the argument.

**3.1. Uniformization Theorem.** *Let $X$ be an affinoid manifold (i.e., all its local rings are regular) and assume $K$ has characteristic zero. Let $\Sigma$ be a subanalytic subset of $X$. Then there exists a finite collection $E$ of maps $\pi \colon X_\pi \to X$, with each $X_\pi$ again affinoid, such that the following properties hold.*

(i) *Each $\pi \in E$ is the composition $\psi_1 \circ \cdots \circ \psi_m$ of finitely many local blowing up maps $\psi_i$ with nowhere dense and smooth centre, for $i < m$.*
(ii) *The union of all the $\operatorname{Im}(\pi)$, for $\pi \in E$, equals $X$.*
(iii) *For each $\pi \in E$, we have that $\pi^{-1}(\Sigma)$ is semianalytic in $X_\pi$.*

*Proof.* Let $X = \operatorname{Sp} A$. As already mentioned, we will use Embedded Resolution of Singularities on $X$ and more particularly the following corollary to it: given $p, q \in A$, then there exists a finite collection $E'$ of maps, such that (i) and (ii) hold, for each $\pi \colon X_\pi \to X$ in $E'$, and furthermore either $p \circ \pi$ divides $q \circ \pi$, or vice versa, $q \circ \pi$ divides $p \circ \pi$, in the affinoid algebra of $X_\pi$. See for instance [**Sch 2**, Lemma 4.2] for a proof.

From our Quantifier Elimination (**2.7**), we know that $\Sigma$ is **D**-semianalytic. By a (not too difficult) argument, involving an induction on the number of times the function **D** appears in one of the describing functions of $\Sigma$ (for details see [**Sch 2**, Theorem 4.4] ), we can reduce to the case that there is only one such occurrence. In other words, we may assume that there exist a quantifier free formula $\psi(\bar{\boldsymbol{x}}, \boldsymbol{y})$ in the language $\mathcal{L}_{\text{an}}$ and functions $p, q \in A$, such that $x \in \Sigma$, if and only if,

$$\psi(x, \mathbf{D}(p(x), q(x))) \quad \text{holds.} \tag{1}$$

After an appeal to the aforementioned corollary of Embedded Resolution of Singularities to $p$ and $q$, and since we only seek to prove our result modulo finite collections of maps for which (i) and (ii) holds, we may already assume that either $p$ divides $q$ or $q$ divides $p$. In the former case, there is some $h \in A$, such that $q = hp$ in $A$. Therefore, $\mathbf{D}(p(x), q(x)) = 0$, unless $q(x) \neq 0$ and $|h(x)| = 1$, in which case it is equal to $1/h(x)$. Let $U_1$ be the affinoid subdomain defined by $|h(x)| \leq 1/2$ and $U_2$ by $|h(x)| \geq 1/2$, so that $\{U_1, U_2\}$ is an admissible affinoid covering of $X$. Hence $x \in U_1$ belongs to $\Sigma$, if and only if, $\psi(x, 0)$ holds and $x \in U_2$ belongs to $\Sigma$, if and only if,

$$[|h(x)| \geq 1 \wedge q(x) \neq 0 \wedge \psi(x, 1/h(x))] \vee [(|h(x)| < 1 \vee q(x) = 0) \wedge \psi(x, 0)]$$

holds. Observe that $1/h$ belongs to the affinoid algebra of $U_2$, since $h$ does not vanish on $U_2$. In other words, $\Sigma$ is semianalytic on both sets and whence on the whole of $X$.

In the remaining case that $q$ divides $p$, i.e. there is some $h \in A$, such that $qh = p$ in $A$, we have an even simpler description of $\Sigma$, namely $x \in \Sigma$, if and only if,

$$[p(x) \neq 0 \wedge \psi(x, h(x))] \vee [p(x) = 0 \wedge \psi(x, 0)]$$

holds, again showing that $\Sigma$ is semianalytic. ∎

**3.2. Corollary.** *Suppose $K$ has characteristic zero and let $\Sigma \subset R^2$. If $\Sigma$ is subanalytic, then in fact it is semianalytic.*

*Proof.* In [**Sch 3**, Theorem 3.2] this is proved for the subclass of strongly subanalytic sets. However, in its proof, nowhere we have made essential use of the strongness (=overconvergency), and hence the same proof applies verbatim (see also the final remark in the introduction of loc. cit.). ∎

*Remark.* Using ABHYANKAR's Embedded Resolution of Singularities in positive characteristic for excellent local rings of dimension two, one can remove the assumption on the characteristic in the above Corollary.

## 4. Elimination along Flat Maps

**4.1. Definition.** This section will be devoted to a proof of (**2.2**). In it, we will need some properties of the *reduction functor* applied to an affinoid algebra. However, for our purposes, we do not need to introduce the whole machinery of reductions but can do with an ad hoc construction to be presented below. First of all, let's fix some notation. As before, $R$ denotes the valuation ring of $K$, i.e., all $r \in K$ with $|r| \leq 1$, and $\wp$ will denote the maximal ideal of $R$, i.e., all $r \in K$, such that $|r| < 1$. The residue field $R/\wp$ will be denoted by $\bar{K}$. Notice that it is also an algebraically closed field.

We will call an $R$-algebra $A^\circ$ an *admissible* algebra, if $A^\circ$ is flat as an $R$-algebra and *topologically of finite type*, meaning of the form $R\langle S \rangle / I^\circ$, for some finitely generated ideal $I^\circ$ and some variables $S = (S_1, \ldots, S_N)$. From a given admissible algebra $A^\circ$, we can construct an affinoid algebra by tensoring over $K$, namely let $A = A^\circ \otimes_R K$. Flatness now guarantees that $A^\circ \subset A$.

If we start with an affinoid algebra $A = K\langle S \rangle / I$ and define $A^\circ$ as $R\langle S \rangle / I^\circ$ with $I^\circ = I \cap R\langle S \rangle$ as above, then $I^\circ$ is finitely generated and $A^\circ$ is torsion-free whence flat over $R$, that is to say, $A^\circ$ is admissible. By tensoring over $K$ we recover our original affinoid algebra, i.e., $A = A^\circ \otimes_R K$. However, $A^\circ$ depends on the particular choice of representing $A$ as a quotient of some $K\langle S \rangle$.

For the sake of simplicity, let us assume that $K$ is algebraically closed (this assumption is not essential, although the proofs would require some modifications for the general case). Let $A^\circ$ be an admissible $R$-algebra and let $A = A^\circ \otimes_R K$ be the corresponding affinoid algebra. With respect to the structure map $R \to A^\circ$, any prime ideal of $A^\circ$ lies either above $(0)$ or above $\wp$. The former prime ideals are in one-one correspondence with the prime ideals of $A$. Hence, in particular, we will consider $\operatorname{Sp} A$ as a subset of $\operatorname{Spec}(A^\circ)$ (via the canonical map $\operatorname{Sp} A \to \operatorname{Spec}(A^\circ)$ induced by $A^\circ \subset A$). Let us call a map $x^\circ \colon \operatorname{Spec} R \to \operatorname{Spec}(A^\circ)$ an *R-rational point*. This means that the image of the generic point of $\operatorname{Spec} R$ (i.e., the prime ideal $(0)$) is a prime ideal $\mathfrak{p}^\circ$ of $A^\circ$, such that $A^\circ/\mathfrak{p}^\circ \cong R$. One easily verifies that $\mathfrak{p}^\circ$, or, rather, $\mathfrak{p}^\circ A$, lies in $\operatorname{Sp} A$ (under the above identification), and hence induces a $K$-rational point $x \colon \operatorname{Spec} K \to \operatorname{Sp} A$ via the map $A/\mathfrak{p}^\circ A \cong K$. (As customary, we will further identify $x$ with the unique point of $\operatorname{Sp} A$ it determines). Conversely, any maximal ideal $\mathfrak{m}$ of $A$ (or, equivalently, any point $x \in \operatorname{Sp} A$) determines an $R$-rational point $x^\circ$ given by the map $A^\circ/(\mathfrak{m} \cap A^\circ) \cong R$. In other words, to give a point $x \in \operatorname{Sp} A$ is the same as to give an $R$-rational point $x^\circ$.

Now, given a point $x \in \operatorname{Sp} A$, let $x^\circ$ be the corresponding $R$-rational point and denote by $\bar{x}$ the map $\operatorname{Spec} \bar{K} \to \operatorname{Spec}(A^\circ)$ obtained by restricting $x^\circ$ to the closed immersion $\operatorname{Spec} \bar{K} \hookrightarrow \operatorname{Spec} R$. We call $\bar{x}$ the *reduction* of $x$ (or, $x^\circ$). In other words, if $\mathfrak{p}^\circ$ is the prime ideal of $A^\circ$ given as the image of the generic point under $x^\circ$, then $\bar{x}$ is determined by the prime ideal $\mathfrak{p}^\circ + \wp A^\circ$. This is in fact a maximal ideal of $A^\circ$. Let us denote the maximal spectrum of $A^\circ$ by $\operatorname{Max}(A^\circ)$. The *reduction map* $\xi \colon \operatorname{Sp} A \to \operatorname{Max}(A^\circ)$ is the map given by sending $x$ to its reduction $\bar{x}$. A word of caution: the reduction map is not induced by any algebra morphism.

If we rather view $\operatorname{Sp} A$ as a subset of $R^N$ than as a maximal spectrum, we can identify $\xi(x)$ with the maximal ideal of $A^\circ$ consisting of all $p \in A^\circ$ for which $|p(x)| < 1$. The reduction map $\xi$ is functorial in the following sense. Let $\varphi^\circ \colon A^\circ \to B^\circ$ be an $R$-algebra morphism of admissible algebras and let $\varphi \colon A \to B$ be the morphism of affinoid algebras obtained by tensoring $\varphi^\circ$ with $K$. Then we have a commutative diagram

$$(\diamond) \quad \begin{array}{ccc} \operatorname{Sp} B & \xrightarrow{f} & \operatorname{Sp} A \\ \xi \downarrow & & \downarrow \xi \\ \operatorname{Max}(B^\circ) & \xrightarrow{f^\circ} & \operatorname{Max}(A^\circ) \end{array}$$

where $f$ and $f^\circ$ are the respective maps on the maximal spectra induced by $\varphi$ and $\varphi^\circ$.

It is well-known (see for instance [**Meh**]) that the reduction map is surjective (regardless whether $K$ is algebraically closed or not). For the sake of the reader's convenience we have added a proof of this fact. It is an immediate consequence of the lemma below, which we need anyway for our proof of RAYNAUD's result. In it, the relevance of flatness is made apparent.

**4.2. Lemma.** *Let $A^\circ \to B^\circ$ be a flat $R$-algebra morphism of admissible $R$-algebras and let $f \colon Y^\circ = \operatorname{Spec}(B^\circ) \to X^\circ = \operatorname{Spec}(A^\circ)$ denote the corresponding map of affine schemes. Let $x^\circ \colon \operatorname{Spec} R \to X^\circ$, be an $R$-rational point of $X^\circ$ and let $\bar{x}$ denote its reduction $\operatorname{Spec} \bar{K} \to X^\circ$. Suppose there exists a $\bar{K}$-rational point $\bar{y} \colon \operatorname{Spec} \bar{K} \to Y^\circ$, such that*

$$(1) \quad \begin{array}{ccc} \operatorname{Spec} \bar{K} & \xrightarrow{\bar{y}} & Y^\circ \\ \| & & \downarrow f \\ \operatorname{Spec} \bar{K} & \xrightarrow{\bar{x}} & X^\circ \end{array}$$

*commutes. Then there exists an $R$-rational point $y^\circ$ of $Y^\circ$, which has reduction $\bar{y}$, and is such that*

$$(2) \quad \begin{array}{ccc} \operatorname{Spec} R & \xrightarrow{y^\circ} & Y^\circ \\ \| & & \downarrow f \\ \operatorname{Spec} R & \xrightarrow{x^\circ} & X^\circ \end{array}$$

commutes. We call $y^\circ$ a factorisation of $x^\circ$ lifting $\bar{y}$.

*Proof.* Let $\mathfrak{p}^\circ$ be the prime ideal of $A^\circ$ associated to $x^\circ$ (i.e., the image of the generic point under $x^\circ$). Let $\bar{\mathfrak{p}} = \mathfrak{p}^\circ + \wp A^\circ$, so that it is the maximal ideal of $A^\circ$ associated to $\bar{x}$. Finally, let $\bar{\mathfrak{q}}$ be the maximal ideal of $B^\circ$ associated to $\bar{y}$, so that the commutativity of (1) translates into

$$\bar{\mathfrak{p}} = \bar{\mathfrak{q}} \cap A^\circ. \tag{3}$$

Since $A^\circ \to B^\circ$ is flat, the Going Down Theorem (see for instance [**Mats**, Theorem 9.5]) guarantees the existence of a prime ideal $\mathfrak{n}^\circ$ of $B^\circ$, such that $\mathfrak{n}^\circ \subset \bar{\mathfrak{q}}$ and

$$\mathfrak{p}^\circ = \mathfrak{n}^\circ \cap A^\circ. \tag{4}$$

Let $\mathfrak{q}^\circ$ be an ideal of $B^\circ$, maximal with respect to the following two conditions

$$\mathfrak{n}^\circ \subset \mathfrak{q}^\circ \subset \bar{\mathfrak{q}} \tag{5}$$
$$\mathfrak{q}^\circ \cap R = (0). \tag{6}$$

The reader easily verifies that such an ideal is necessarily a prime ideal. Moreover, by (3), (4) and (5), we must have inclusions

$$\mathfrak{p}^\circ \subset \mathfrak{q}^\circ \cap A^\circ \subset \bar{\mathfrak{p}}.$$

In view of (6), the latter of these must be strict. Since $A^\circ/\mathfrak{p}^\circ \cong R$, the only two prime ideals of $A^\circ$ containing $\mathfrak{p}^\circ$ are $\mathfrak{p}^\circ$ itself and $\bar{\mathfrak{p}}$. Therefore, we conclude that

$$\mathfrak{q}^\circ \cap A^\circ = \mathfrak{p}^\circ. \tag{7}$$

Let $B = B^\circ \otimes_R K$ be the associated affinoid algebra. We claim that $\mathfrak{q}^\circ B$ is a maximal ideal of $B$. Assuming the claim, we have an inclusion of $R$-algebras

$$R \hookrightarrow B^\circ/\mathfrak{q}^\circ \hookrightarrow B/\mathfrak{q}^\circ B \cong K.$$

Again the last inclusion must be strict and since $R$ is a rank-one valuation ring, the first inclusion is in fact an isomorphism. In other words, if $y^\circ$ is the point of $Y^\circ$ corresponding to $\mathfrak{q}^\circ$, then it is an $R$-rational point. Moreover, $\mathfrak{q}^\circ + \wp B^\circ$ is then a maximal ideal, containing $\mathfrak{q}^\circ$ and contained in $\bar{\mathfrak{q}}$ by (5), and hence equal to the latter. This shows that $y^\circ$ has $\bar{y}$ as its reduction, as required.

It remains to prove the claim. To this end, let $b \in B$ not belonging to $\mathfrak{q}^\circ B$. We can find $0 \neq \pi \in \wp$, such that $\pi b$ belongs to $B^\circ$ and even to $\wp B^\circ$. In particular, it belongs to $\bar{\mathfrak{q}}$. By the maximality of $\mathfrak{q}^\circ$, we must have that

$$(\mathfrak{q}^\circ + \pi b B^\circ) \cap R \neq (0).$$

From this it follows that

$$\mathfrak{q}^\circ B + bB = (1),$$

showing that $B/\mathfrak{q}^\circ B$ is a field, as wanted. ∎

**4.3. Corollary.** *Let $B^\circ$ be an admissible algebra and $B = B^\circ \otimes_R K$. Then the map $\xi \colon \mathrm{Sp}\, B \to \mathrm{Max}(B^\circ)$ is surjective.*

*Proof.* Let $\bar{\mathfrak{m}}$ be a maximal ideal of $B^\circ$ and let $\bar{y}$ be the corresponding $\bar{K}$-rational point of $Y^\circ = \mathrm{Spec}(B^\circ)$. Apply the previous lemma to $A = R$ and $x^\circ \colon \mathrm{Spec}\, R \to \mathrm{Spec}\, R$ given by the identity morphism on $R$. Note that by assumption, $R \to B^\circ$ is indeed flat so that the lemma applies. Let $y^\circ$ be a factorization of $x^\circ$ lifting (i.e., with reduction) $\bar{y}$. This means precisely that the corresponding point $y \in \mathrm{Sp}\, B$ has reduction $\xi(y) = \bar{y}$, as required. ∎

**4.4. Corollary.** *Let $A^\circ$ be an admissible algebra with corresponding affinoid algebra $A = A^\circ \otimes_R K$. Let $h_1, \ldots, h_s \in A^\circ$. Let $\Sigma$ denote the semianalytic set of all $y \in \operatorname{Sp} A$, such that $|h_i(y)| < 1$, for $i < r$ and $|h_i(y)| \geq 1$, for $r \leq i < s$. Call such a set* special. *Let $\Sigma^\circ$ denote the locally closed subset of $\operatorname{Max}(A^\circ)$ consisting of all maximal ideals $\bar{\mathfrak{m}}$, such that $h_i \in \bar{\mathfrak{m}}$, for $i < r$, and $h_i \notin \bar{\mathfrak{m}}$, for $r \leq i < s$. Then these two sets are related to one other by*

$$\xi^{-1}(\Sigma^\circ) = \Sigma \qquad \text{and} \qquad \Sigma^\circ = \xi(\Sigma). \tag{1}$$

*In other words, $\xi$ induces a bijection between the class of finite Boolean combinations of special subsets of $\operatorname{Sp} A$ and the class of locally closed subsets of $\operatorname{Spec} \bar{A}$, where $\bar{A} = A^\circ/\wp A^\circ$.*

*Proof.* This is only a matter of writing out the definitions using the surjectivity of $\xi$. ∎

**4.5. Proof of Theorem (2.2).** So we are given a flat map $f \colon \operatorname{Sp} B \to \operatorname{Sp} A$ of affinoid varieties and a special set $\Sigma$ of $\operatorname{Sp} B$, given as

$$\{\, y \in \operatorname{Sp} B \mid |h_i(y)| \diamond_i 1 \text{ for } i < s \,\}, \tag{1}$$

where $h_i \in B$ are of supremum norm at most one and $\diamond_i$ is either $<$ or $\geq$. We want to prove that $f(\Sigma)$ is semianalytic.

Using **(4.6)** below, we may reduce to the case that there exist admissible algebras $A^\circ$ and $B^\circ$, such that $A = A^\circ \otimes_R K$ and $B = B^\circ \otimes_R K$, and there exists a flat morphism of $R$-algebras $A^\circ \to B^\circ$ which induces the map $f$ (after tensoring with $K$), such that $h_i \in B^\circ$, for all $i < t$. Now, let $\Sigma^\circ$ be the locally closed set of $\operatorname{Max}(B^\circ)$ defined by

$$\Sigma^\circ = \{\, \mathfrak{n}^\circ \in \operatorname{Max}(B^\circ) \mid h_i \bowtie_i \mathfrak{n}^\circ \text{ for } i < s \,\},$$

where $\bowtie_i$ stands for $\in$ when $\diamond_i$ is $<$, and $\bowtie_i$ stands for $\notin$ when $\diamond_i$ is $\geq$. By **(4.4)**, we have that $\xi^{-1}(\Sigma^\circ) = \Sigma$. We also observed above that $\operatorname{Max}(B^\circ)$ can be identified with $\operatorname{Max}(\bar{B})$, where $\bar{B} = B^\circ/\wp B^\circ$ and hence we can consider $\Sigma^\circ$ as a subset of the latter space as well. If we also put $\bar{A} = A^\circ/\wp A^\circ$, then since both rings are finitely generated $\bar{K}$-algebras, we can invoke CHEVALLEY's Theorem to conclude that the image of $\Sigma^\circ$ under the induced map $\bar{f} \colon \operatorname{Max}(\bar{B}) \to \operatorname{Max}(\bar{A})$ is a constructible set $\Omega^\circ$. Identifying $\operatorname{Max}(A^\circ)$ with $\operatorname{Max}(\bar{A})$, we may consider $\Omega^\circ$ as a constructible set of the former space as well and as such it is the image of $\Sigma^\circ$ under the map $f^\circ$ induced by $f$. Let $\Omega = \xi^{-1}(\Omega^\circ)$, then using **(4.4)** once more, we obtain that $\Omega$ is semianalytic in $\operatorname{Sp} A$. Hence we proved our theorem once we showed that

$$f(\Sigma) = \Omega. \tag{2}$$

The commutative diagram ($\diamond$) of **(4.1)** expressing the functoriality of $\xi$, provides the inclusion $f(\Sigma) \subset \Omega$, so we only need to deal with the opposite inclusion.

Hence let $x \in \Omega$. Let $x^\circ$ be the corresponding $R$-rational point and let $\bar{x}$ be the reduction $\xi(x)$ of $x$. By assumption, $\bar{x} \in \Omega^\circ$ and hence it is the image under $\bar{f}^\circ$ of some point $\bar{y} \in \Sigma^\circ$. We can apply **(4.2)** to this situation to obtain an $R$-rational point $y^\circ$ lifting $x^\circ$ and compatible with $\bar{y}$. In other words, if $y \in \operatorname{Sp} B$ denotes the point corresponding to $y^\circ$, then this translates into $f(y) = x$ and $\xi(y) = \bar{y}$. In view of **(4.4)** and the fact that $\bar{y} \in \Sigma^\circ$, the latter implies that $y \in \Sigma$, as required. ∎

**4.6. Proposition.** *Let $f \colon Y = \operatorname{Sp} B \to X = \operatorname{Sp} A$ be a flat map of affinoid varieties and let $h_j \in B$, for $j < t$, be of supremum norm at most one. There exist finite coverings $\{U_i = \operatorname{Sp} A_i\}_{i<s}$ of $X$ and $\{V_i = \operatorname{Sp} B_i\}_{i<s}$ of $Y$ by rational subdomains and $R$-algebra morphisms $\varphi_i^\circ \colon A_i^\circ \to B_i^\circ$ of admissible algebras, such that, for all $i < s$, we have that*

(1) $A_i = A_i^\circ \otimes_R K$ and $B_i = B_i^\circ \otimes_R K$,
(2) *the morphism $\varphi_i^\circ$ induces the map $f|_{V_i} \colon V_i \to U_i$,*
(3) *the morphism $\varphi_i^\circ$ is flat,*
(4) $h_j \in B_i^\circ$, *for all $j < t$.*

*Proof.* Since the $h_j$ are of norm at most one and using [**BGR**, 6.4.3. Theorem 1], we can find an admissible algebra $B^\circ$ containing all $h_j$ with $B^\circ \otimes_R K = B$, an admissible algebra $A^\circ$ with $A^\circ \otimes_R K = A$ and an $R$-algebra morphism $\varphi^\circ \colon A^\circ \to B^\circ$ inducing the map $f$. In general, $\varphi^\circ$ will not be flat. To remedy this, we use [**Meh**, 3.4.8], in order to find admissible coverings as asserted, for which (1)-(3) holds. Moreover, from the proof in loc. cit., it follows that $B_i^\circ$ is a quotient of $A_i^\circ \otimes_{A^\circ} B^\circ$. Therefore also (4) is satisfied. ∎

*Remark.* The result in MEHLMANN's paper is quite an intricate matter, using RAYNAUD's approach on rigid analysis through formal schemes and admissible formal blowing ups; an alternative proof can be found in [**BL**].

## Appendix

**Osgood's Counterexample.**

Suppose $K$ has characteristic zero. Let $\Sigma$ be the subset of $R^3$ given by the triples $(s, st, se^{\tau t})$, where $s, t \in R$. Here $\tau \in R$ with $0 < |\tau| < 1$, so that $e^{\tau T} = \sum_i (\tau T)^i/i!$ is indeed a strictly convergent power series over $K$. We will show the following two properties of $\Sigma$.

**Facts.**
(i) *The set $\Sigma$ is subanalytic but not semianalytic as a subset of $R^3$.*
(ii) *The set $\Sigma \setminus \{(0,0,0)\}$ is semianalytic as a subset of $R^3 \setminus \{(0,0,0)\}$.*

This proves that (i) not every subanalytic set is semianalytic, i.e., $R$ does not admit elimination of quantifiers in the language $\mathcal{L}_{\text{an}}$; and (ii) being semianalytic or subanalytic is a property of how a set embeds in a larger set, rather than an intrinsic property of that set, since $\Sigma \setminus \{(0,0,0)\}$ is semianalytic in $R^3 \setminus \{(0,0,0)\}$ but not in $R^3$. Moreover, from the proof it also follows that (vi) in **(2.6)** is false for general subanalytic sets.

Let us first show (i). By (v) of **(2.6)**, we must have that $\Sigma$ has dimension at most two, since $\Sigma$ is equal to the image of the map

$$\gamma\colon R^2 \to R^3\colon (s,t) \mapsto (s, st, se^{\tau t}).$$

However, if $G(U_1, U_2, U_3) = \sum_i G_i(U_1, U_2, U_3) \in K\langle U_1, U_2, U_3 \rangle$ is a power series vanishing on $\Sigma$, where each $G_i$ is homogeneous of degree $i$, then $\sum_i s^i G_i(1, t, e^{\tau t}) = 0$, for all $s, t \in R$, whence all $G_i(1, t, e^{\tau t}) = 0$, for all $i$, which implies that each $G_i = 0$, since the three functions $1, T$ and $e^{\tau T}$ are algebraically independent. In other words, the Zariski closure of $\Sigma$ equals $R^3$, so that by (vi) of loc. cit., if $\Sigma$ were semianalytic, its dimension would be 3, contradiction.

Let us now prove (ii). If $W_n$ denotes the admissible affinoid of $R^3 \setminus \{(0,0,0)\}$ given by $|u_i| \geq |2|^{-n}$, for $i = 1, 2, 3$, then the collection of all such $W_n$, with $n = 0, 1, \ldots$, gives an admissible affinoid covering of $R^3 \setminus \{(0,0,0)\}$. Hence to prove (ii), it will be sufficient to prove that $\Sigma \cap W_n$ is semianalytic in $W_n$. Let $\tau' \in R$ be such that $1 > |\tau'| > |\tau|$. Let $W_{n,1}$ be the affinoid subdomain of $W_n$ given by $|u_2 \tau'| \leq |u_1|$ and $W_{n,2}$ the subdomain given by $|u_2 \tau'| \geq |u_1|$, so that $\{W_{n,1}, W_{n,2}\}$ is an admissible covering of $W_n$. If $(u_1, u_2, u_3) \in \Sigma$, then in particular $|u_2| \leq |u_1|$, so that $\Sigma \cap W_{n,2} = \emptyset$. Hence it only remains to verify the semianalyticity of $\Sigma$ when restricted to $W_{n,1}$. But then $(u_1, u_2, u_3) \in W_{n,1}$ belongs to $\Sigma$, if and only if,

$$|u_2| \leq |u_1| \wedge u_3 = e^{\tau u_2/u_1} u_1.$$

This is indeed a semianalytic description, in view of the fact that $e^{\tau u_2/u_1}$ belongs to the affinoid algebra of $W_{n,1}$ since $\tau u_2/u_1 = (\tau/\tau')(\tau' u_2/u_1)$ and $|\tau/\tau'| < 1$.

One could as easily give an example as in (i) in positive characteristic, since the only property we used of the series $e^{\tau T}$ is that it is algebraically independent from 1 and $T$. Other instances of (ii) are given by any situation where $\pi\colon \tilde{X} \to X$ is a finite sequence of local blowing up maps and $\Sigma$ is subanalytic (non-semianalytic) such that $\pi^{-1}(\Sigma)$ is semianalytic, since outside a nowhere dense closed analytic subset, both sets are isomorphic, as already observed in the remark following **(2.3)**. (Also the remark after **(2.7)** is relevant here). However, the example in (ii) where a set becomes semianalytic after taking away a single point out of the space, we seem to have made essential use of the fact that $e^{\tau T}$ is overconvergent and we do not know of any such example not using overconvergency.

MATHEMATICAL INSTITUTE
UNIVERSITY OF OXFORD
24-29 ST. GILES
OXFORD OX1 3LB (UNITED KINGDOM)
*E-mail address*: gardener@maths.ox.ac.uk

FIELDS INSTITUTE
222 COLLEGE STREET
TORONTO, ONTARIO, M5T 3J1 (CANADA)
*E-mail address*: hschoute@fields.utoronto.ca